\documentclass[11pt]{amsart}
\usepackage{amsmath, amssymb,amsthm,enumerate}
\usepackage{hyperref}
\newcommand{\cal}{\mathcal}

\newtheorem{thm}{Theorem}
\newtheorem{prop}[thm]{Proposition}
\newtheorem{lema}[thm]{Lemma}
\newtheorem{defi}[thm]{Definition}
\newtheorem{cor}[thm]{Corollary}

\newcommand{\bY}{\boldsymbol{Y}}

\title[Graphical independence structure in DPP's]{Graphical structure of 
conditional independencies in determinantal 
point processes}
\author{Tvrtko Tadi\'c}

\address{Department of Mathematics, University of Washington, Seattle \and 
Department of Mathematics, University of Zagreb, Zagreb, Croatia}
\email{tadic@math.washington.edu}
\urladdr{www.math.washington.edu/$\sim$tadic}
\email{tvrtko@math.hr}
\begin{document}

 \begin{abstract}
  Determinantal point process have recently been used as models in machine learning and
this has raised questions regarding the characterizations of conditional independence.
In this paper we investigate characterizations of conditional independence. 
We describe some conditional independencies
through the conditions on the kernel of a determinantal point process, and show many can 
be obtained using the graph induced by a kernel of the $L$-ensemble.
 \end{abstract}
\maketitle
In recent years there have been several machine learning papers
about the applications of determinantal point processes (DPP's) \cite{doc_sort},
\cite{map}, \cite{l_dpp}, \cite{k-dpp}\ldots An overview of theory, recent applications and problems 
in learning DPP's is given in a recent extensive survey \cite{kul_task} by Kulesza and Taskar.\vspace{0.2cm} 

In a private communication with Ben Taskar, one of the questions from survey \cite{kul_task} (see \S 7.3),  that remains for future research, was brought to author's attention:
\begin{itemize}
 \item {\it Is there a simple characterization of the conditional independence relations encoded by a DPP?}
\end{itemize}

This question arises naturally having in mind conditional 
independence structure models (see \cite{indp_struct}), such as graphical models (see \cite{grph:model})
that are often used. 
\vspace{0.2cm}

It turns out that, from the mathematical view point, elegant  characterizations, similar to those in
graphical models, exist. 
This paper provides two (main) characterizations:
\begin{itemize}
 \item the block in a Schur complement of the kernel has to be a $0$-block (Theorem \ref{thm:mainstrong}, Proposition \ref{prop:main});
 \item we can use the structure of the graph induced by the kernel of the $L$-ensemble 
to read many conditional independencies in the process (Theorem \ref{thm:Lgrph}, Proposition \ref{prop:Lgrph}).
\end{itemize}

\section{Introduction to the model}

In this paper $K$ will be a positive semi-definite $N\times N$ matrix. Let ${\bf 0}\preceq K\preceq {\bf I}$,
 ${\cal Y}=\{1,\ldots, N\}$. We call a random subset  $\boldsymbol{Y}$ 
of ${\cal Y}$ a \textbf{determinantal point process} if the following holds
$$\Pr(A\subset \boldsymbol{Y})=\det(K_A),$$
and by definition $\Pr(\emptyset \subset \boldsymbol{Y})=1$. (Where $K_{A}=[K_{ij}]_{i,j\in A}$.)\vspace{0.2cm}

Basically, we have a set of $N$ points, and we pick a random subset $\bY$ of them.
We model the probability that all the points in the subset $A$ were chosen by 
$\det( K_A)$. \vspace{0.2cm}

Instead of modeling with the \textbf{kernel} $K$, in practice a determinantal point process is modeled as an $L$-ensemble.
The process $\bY$ is called the \textbf{$L$-ensemble} with the \textbf{kernel} $L$ if 
$$\Pr(\boldsymbol{Y}=A)=\frac{\det(L_A)}{\det(L+I)},$$
where $L$ is a positive semi-definite matrix.

\thm{An $L$-ensemble with kernel $L$ is a DPP with the kernel 
$$K=L(L+I)^{-1}= I- (L+I)^{-1}.$$ }

\cor{For ${\bf 0}\prec K \prec {\bf I}$, a DPP with a kernel $K$ is an $L$-ensemble where 
\begin{equation}
L=K(I-K)^{-1}=(I-K)^{-1}-I. \label{eq:L} 
\end{equation}

}\rm

The following proposition summarizes some useful results about DPP's (they are all 
proven in \cite{kul_task}). Through this text $K_{AB}=[K_{ij}:i\in A,j\in B]$.

\prop{\label{prop:sumrez}Let $\bY$ be a DPP over ${\cal Y}$ with kernel $K$ and $A\subset {\cal Y}$.
\begin{enumerate}[(a)]
 \item The process $\bY_A=\bY\cap A$ is a DPP with kernel $K_A$.
 \item We have 
$$\Pr(A\subset \bY,B\cap\bY=\emptyset )=(-1)^{|B|}\det\left[\begin{array}{cc}
                       K_A & K_{AB}\\
		      K_{AB}^T & K_B-I
                      \end{array}\right].$$
 \item The process ${\cal Y}\setminus \bY$ is a DPP with the kernel $I-K$.
\end{enumerate}}\rm

\vspace{0.1cm}

For more on results and properties of DPP's see \cite{borodin} or \S4 in \cite{dtpr_sv}.\vspace{0.1cm}

In further text, we will assume ${\bf 0}\prec K \prec {\bf I}$ and ${\bf 0}\prec L$.

\section{Independencies}

Under which conditions for three
disjoint subsets $A,B,C$ of ${\cal Y}$ we have\footnote{We use the notation $S_1\perp S_2|S_3$
to denote that $S_1$ is independent of $S_2$ given $S_3$.} 
\begin{equation}
 (A\subset \boldsymbol{Y})\perp (B\subset \boldsymbol{Y})\, | \, (C\subset \boldsymbol{Y}).\label{in1}
\end{equation}
This was investigated by Kulesza in \cite{kul}, where the answer is given
for the case $|A|=|B|=1$. We will give a very general answer in Proposition \ref{prop:main}.

\subsection{Independence in DPPs}
We will start with the case $C=\emptyset$. When is
\begin{equation}
 (A\subset \boldsymbol{Y})\perp (B\subset \boldsymbol{Y})?\label{in2}
\end{equation}

The following are some known  technical results from matrix analysis (see \cite{horn_john}).

\lema{\label{lem:2}Let 
\begin{equation}
 M_+=\left[\begin{array}{cc}
     U & V\\
     {V}^T & W\\
    \end{array}\right]\quad and \quad M_-=\left[\begin{array}{cc}
     U & -V\\
     -{V}^T & W\\
    \end{array}\right] \label{def:m+}
\end{equation}
be quadratic matrices.
\begin{enumerate}[(a)]
       \item If $M_+$ and $M_-$
are symmetric matrices their eigenvalues are the same with the same multiplicity.
Further their determinants are also the same.
	\item $M_+$ is positive definite if and only if $M_-$ is positive definite.
        \item $M_+$ is positive definite, if and only if
\begin{equation}
	  U-VW^{-1}V^T \quad {\rm and}\quad W \label{psdschrchar}
\end{equation}
are positive definite.

\item \label{lem:shrc} If $W$ is non-singular, then 
$$\det (M_+)=\det (M_-)= \det (W) \det (U-VW^{-1}V^T)  $$

      \end{enumerate}}
\cor{\label{cor:0}If $M_+$ is a positive (semi)defnite matrix so is 
  $$M_0=\left[\begin{array}{cc}
     U & {\bf 0}\\
     {\bf 0}^T & W\\
    \end{array}\right].$$
\begin{proof}
 Follows from the fact that $M_0=\frac{1}{2}(M_++M_-)$.
\end{proof}
}

\rm 

                                                    
We following technical lemma will be the key for conditional independencies.
\lema{\label{lem:main}Let $A$ be a positive definite and $B$ a positive semi-definite $N\times N$ matrices.
If $\det(A+B)=\det A$, then $B=0$. 
\begin{proof}
 Since $A$ is positive definite, there exists a positive definite matrix $\sqrt{A}$, such that
$A=(A^{1/2})^2$. Therefore, since $\det A= (\det A^{1/2})^2$, 
we have 
\begin{equation}
\label{eq:det1} \det(I+A^{-1/2}BA^{-1/2})=1 .
\end{equation}

It is not hard to see that 
$A^{-1/2}BA^{-1/2}$ is a positive semi-definite matrix.
Hence $(\ref{eq:det1})$ is equivalent (using the eigenvalue decomposition)
$$(1+\lambda_1)\ldots (1+\lambda_N)=1,$$
where $\lambda_1,\ldots,\lambda_N$ are eigenvalues of $A^{-1/2}BA^{-1/2}$.
Since this matrix is positive semi-definite, $\lambda_j\geq 0$ for $j=1,\ldots, N$ and therefore 
we have $\lambda_1=\ldots =\lambda_N=0$. Hence, $A^{-1/2}BA^{-1/2}={\bf 0}$ and the claim follows.
\end{proof}
}

\cor{\label{cor:2}Let 
$$M=\left[\begin{array}{cc}
     U & V\\
     {V}^T & W\\
    \end{array}\right].$$
If one of the following conditions holds
\begin{enumerate}[(a)]
 \item $M$ is positive definite;
 \item $U$ is positive definite and $W$ is negative definite;
 \item $M$ is negative definite;
 \item $U$ is negative definite and $W$ is positive definite;
\end{enumerate}
then the equality
$$\det\left[\begin{array}{cc}
     U & V\\
     {V}^T & W\\
    \end{array}\right]= \det U \det W.$$
holds if and only if $V={\bf 0}$.
\begin{proof}
If $V={\bf 0}$ the claim is clear.

We will prove cases (a) and (b), cases (c) and (d) follow from them.

Assume $\det W=\det U\det W$. Using Lemma \ref{lem:2}(\ref{lem:shrc}) we get 
$\det (M)= \det (W) \det (U-VW^{-1}V^T)$, and from the assumption 
we have 
$$\det U=\det (U-VW^{-1}V^T).$$

(a): $A=U-VW^{-1}V^T$ is positive definite (by $(\ref{psdschrchar})$) and $B=VW^{-1}V^T$
is positive semi-definite. By Lemma \ref{lem:main} we have $B={\bf 0}$.
Now, let $V^T=[v_1,\ldots, v_m]$. Since, $B={\bf 0}$, $v_j^TW^{-1}v_j=0$,
and since $W^{-1}$ is positive definite we have $v_j={\bf 0}$ for $j=1\ldots m$.\vspace{0.2cm}

(b): Set $A=U$ and $B=V(-W^{-1})V^T$. Since $-W^{-1}$ is positive definite,
$B$ is positive semi-definite  and by Lemma \ref{lem:main}, $B={\bf 0}$. 
Using the same approach as in case (a) we get $V={\bf 0}$.
\end{proof}
}

\thm{\label{thm:ab}If $K$ is a kernel for the determinantal point process
$\bY$ over ${\cal Y}$, $A$ and $B$ disjoint subsets of ${\cal Y}$, then 
$(A\subset \bY)\perp (B\subset \bY)$ if and only if $K_{AB}={\bf 0}$.
\begin{proof}
 By definition, we have $(A\subset \bY)\perp (B\subset \bY)$ if and only if
$$\Pr(A\cup B\subset \bY)=\Pr((A\subset \bY)\cap (B\subset \bY))=\Pr(A\subset \bY)\Pr(B\subset \bY).$$
This is equivalent to 
$$\det K_{A\cup B}=\det\left[\begin{array}{cc}
                       K_A & K_{AB}\\
		      K_{AB}^T & K_B
                      \end{array}
 \right]= \det K_A \det K_B.$$
By Corollary \ref{cor:2}, this holds if and only if $K_{AB}={\bf 0}$.
\end{proof}
}\rm

\cor{If $K$ is a kernel for the determinantal point process
$\bY$ over ${\cal Y}$, $A$ and $B$ disjoint subsets of ${\cal Y}$, then 
$(A\cap \bY=\emptyset)\perp (B\cap \bY=\emptyset)$ if and only if $K_{AB}={\bf 0}$.
\begin{proof}
By Proposition \ref{prop:sumrez}(c) ${\cal Y}\setminus \bY$ is DPP with kernel $I-K$. 
The claim now follows from Theorem \ref{thm:ab}. 
\end{proof}
}\rm

\thm{\label{thm:ab^c}
If $K$ is a kernel for the determinantal point process
$\bY$ over ${\cal Y}$, $A$ and $B$ disjoint subsets of ${\cal Y}$, then 
$(A\subset \bY)\perp (B\cap \bY=\emptyset)$ if and only if $K_{AB}={\bf 0}$.
\begin{proof}
 By Proposition \ref{prop:sumrez} (b) we know that  $(A\subset \bY)\perp (B\cap \bY=\emptyset)$
if and only if
$$\Pr(A\subset \bY, B\cap \bY=\emptyset)=(-1)^{|B|}\det\left[\begin{array}{cc}
                       K_A & K_{AB}\\
		      K_{AB}^T & K_B-I
                      \end{array}\right]$$
$$=\Pr(A\subset \bY)\Pr(B\cap \bY=\emptyset) = \det K_A (-1)^{|B|}\det(K_B-I).$$
By Corollary \ref{cor:2} this is true if and only if $K_{AB}={\bf 0}$
\end{proof}
}\rm

Using the same techniques as in last proofs, we can prove much more.

\thm{\label{thm:abstrong}If $K$ is a kernel for the determinantal point process
$\bY$ over ${\cal Y}$, $A$ and $B$ disjoint subsets of ${\cal Y}$, then the 
processes 
$\bY_A=\bY\cap A$ and $\bY_B=\bY\cap B$ are independent if and only if $K_{AB}={\bf 0}$.
\begin{proof} If $\bY_A$ and $\bY_B$ are independent, then $(A\subset \bY)\perp (B\subset \bY)$, and hence 
by Theorem \ref{thm:ab} the claim follows.

Let $L^{A\cup B}$ denote the kernel of the $L$-ensemble of the process $\bY\cap (A\cup B)$.
If $K_{AB}={\bf 0}$ we know that for $A_1\subset A$ and $B_1\subset B$ we have
$$\Pr(A\cap \bY=A_1,B\cap \bY=B_1)= \det (L^{A\cup B}_{A_1\cup B_1})=$$$$ =\det L^{A\cup B}_{A_1} \det L^{A\cup B}_{B_1}=\Pr(A\cap \bY=A_1)\Pr(B\cap \bY=B_1),$$
since $$L^{A\cup B}_{A_1\cup B_1}=(I-K_{A_1\cup B_1})^{-1}-I= \left[\begin{array}{cc}
                                                 (I-K_{A_1})^{-1}-I & 0\\
						  0 & (I-K_{B_1})^{-1}-I
                                                 \end{array}
\right].$$
 
\end{proof}
}\rm 

The following proposition summarizes the all the results from this subsection.

\prop{\label{prop:char}For a DPP with the kernel ${\bf 0} \prec K\prec {\bf I}$, and $A$ and $B$
disjoint subsets of ${\cal Y}$ the following statements are equivalent:
\begin{enumerate}[(a)]
 \item $(A\subset \bY)\perp (B\subset \bY)$;
  \item $(A\subset \bY)\perp (B\cap \bY=\emptyset)$;
 \item $(A\cap \bY=\emptyset)\perp (B\cap \bY=\emptyset)$;
 \item $\bY_A\perp \bY_B$;
 \item $K_{AB}=0$.
\end{enumerate}
 }\rm
\noindent\emph{Remark.} One might be tempted to think that if
\begin{equation}
 \label{eq:cnex}(A_1\subset \bY, A_2\cap \bY=\emptyset)\perp (B_1\subset \bY, B_2\cap\bY=\emptyset)
\end{equation}

then $K_{A_1\cup A_2,B_1\cup B_2}={\bf 0}$. However, this doesn't
have to be true. Take $$K=\left[\begin{array}{ccc}
                           0.05 & 0& 0.1\\
			    0	& 0.8 & 0.2\\
			    0.1 &0.2 & 0.6
                          \end{array}\right].
$$
It is not hard to check that ${\bf 0}\prec K\prec {\bf I}$.
Set $A_1=\{1\}$, $A_2=\{2\}$, $C_1=\{3\}$ and $C_2=\emptyset$. Clearly,
$K_{A_1\cup A_2,C_1\cup C_2}\neq {\bf 0}$.
However, by Proposition \ref{prop:sumrez} (b) 
$$\Pr( A_1\subset \bY, A_2\cap \bY=\emptyset, B_1\subset \bY)=-\det\left[\begin{array}{ccc}
                           0.05 & 0& 0.1\\
			    0	& -0.2 & 0.2\\
			    0.1 &0.2 & 0.6
                          \end{array}\right]=0.006,$$
is a product of $\Pr( A_1\subset \bY, A_2\cap \bY=\emptyset)=-\det\left[\begin{array}{ccc}
                           0.05 & 0\\
			    0	& -0.2\\
                          \end{array}\right]=0.01$ and $\Pr(B_1\subset \bY)=0.6$. Hence,
in this case $(\ref{eq:cnex})$ is true.
\subsection{Conditional independencies in DPP's}

It is known that conditioned on the event $(A\subset \bY, B\cap \bY =\emptyset)$ 
the process $\bY$ is a DPP. (See \cite{kul_task} or \cite{borodin}.)

\defi{ If  $M$ is a square matrix and $M_C$ is non-singular then we can define (the Schur complement of $M$)
\begin{equation}
M^C= M_{C^c}-M_{C^c,C}M_C^{-1}M_{C,C^c}=M_{C^c}-M_{C^c,C}M_C^{-1}M_{C^c,C}^T. 
\end{equation}}\rm

\noindent\emph{Remark.} By Lemma \ref{lem:2}(c) if $K$ is positive definite, then $K^C$ is positive definite.
On the other hand, if $K\prec {\bf I}$, then, clearly, ${\bf I}-K^C={\bf I}-K_{C^c}+K_{C^c,C}K_CK_{C,C^c}\succ {\bf 0}$.

\lema{\label{lem:3}For the determinantal point process $\bY$ and  some $C\subset {\cal Y}$ such that $|K_C|>0$, 
for every $A\subset C^c$ we have 
$$\Pr(A\subset \bY| C\subset \bY)= \det (K^C_{A}).$$
Hence $Y\cap C^c$ given $(C\subset \bY)$ is a DPP with the kernel $K^C$.}
\begin{proof}
 By definition,
\begin{align*}
 \Pr(A\subset \bY| C\subset \bY)&=\frac{\Pr(A\subset \bY, C\subset \bY)}{\Pr(C\subset \bY)}=\frac{\Pr(A\cup C\subset \bY)}{\Pr(C\subset \bY)}\\
 &=\frac{\det K_{A\cup C}}{\det K_C}=\frac{1}{\det K_C}\det\left[
\begin{array}{ll}
K_A & K_{AC} \\
K_{AC}^T & K_C\\
    \end{array}
    \right]\\
&\stackrel{\textrm{Lem. \ref{lem:2}(\ref{lem:shrc})}}{=} \det(K_A-K_{AC}K_C^{-1}K_{AC}^T)=\det(K^C_A).
\end{align*}

\end{proof}

\thm{\label{thm:main}For the determinantal point process
$\bY$ over ${\cal Y}$ with the kernel $K$, and $A,B,C$ disjoint subsets of ${\cal Y}$, then
$$(A\subset \boldsymbol{Y})\perp (B\subset \boldsymbol{Y})\, | \, (C\subset \boldsymbol{Y})$$ 
 is true if and only if $K^C_{AB}=0$, i.e.
\begin{equation}
K_{AB} =\begin{cases}
	  K_{AC}K_C^{-1}K_{BC}^T & C\neq \emptyset\\
	   0 & C= \emptyset
        \end{cases}
 \label{eq:2} 
\end{equation}
%
%
\begin{proof}
If $C=\emptyset$ the claim follows from Theorem \ref{thm:ab}.
When $C\neq \emptyset$ from Lemma \ref{lem:3} we know that $\bY \cap C^c|(C\subset \bY)$
is a DPP with kernel $K^C$. Now, by Theorem \ref{thm:ab} $(A\subset \bY)$ and $(B\subset \bY)$
are independent given $(C\subset \bY)$ if and only if $K^C_{AB}={\bf 0}$. Since $K^C_{AB}=K_{AB}-K_{AC}K_C^{-1}K_{CB}$,
the claim follows.\end{proof}
}\rm

Using the same argumentation and Theorem \ref{thm:abstrong} we have the following result.

\thm{\label{thm:mainstrong}If $K$ is a kernel for the determinantal point process
$\bY$ over ${\cal Y}$, and $A,B,C$ disjoint subsets of ${\cal Y}$, then
$\bY_A= \boldsymbol{Y}\cap A$ and $\bY_B=\bY\cap B$ are independent given  $(C\subset \boldsymbol{Y})$ 
 if and only if $K^C_{AB}=0$, i.e. (\ref{eq:2}) is true.}\rm \vspace{0.2cm}

The following is a generalization of the Proposition \ref{prop:char}.

\prop{\label{prop:main}For a DPP with the kernel ${\bf 0} \prec K\prec {\bf I}$, and $A,B,C$
disjoint subsets of ${\cal Y}$ the following statements are equivalent:
\begin{enumerate}[(a)]
 \item $(A\subset \bY)\perp (B\subset \bY)|(C\subset \bY)$;
  \item $(A\subset \bY)\perp (B\cap \bY=\emptyset)|(C\subset \bY)$;
  \item $(A\cap \bY=\emptyset)\perp (B\cap \bY=\emptyset)|(C\subset \bY)$;
 \item $\bY_A\perp \bY_B|(C\subset \bY)$;
 \item $K^C_{AB}=0$.
\end{enumerate}
 }\rm\vspace{0.2cm}

%
It is known (see for example (7.7.5) in \cite{horn_john}) that 
\begin{equation}
 (K^{-1})_{C^c}=(K^C)^{-1}.  \label{eq:com_inv}
\end{equation}

\cor{Let ${\cal Y}$ be a union of disjoint sets $\{i\}$, $\{j\}$ and $C={\cal Y}\setminus\{i,j\}$.
Then $K^{-1}_{ij}=0$ if and only if $K^{C}_{ij}=0$.
\begin{proof}
Note that $K^C$ is a $2\times 2$ matrix. $K^{C}_{ij}=0$ if and only if $K^C$ is 
a diagonal matrix. This is so if and only if $(K^{C})^{-1}_{ij}=0\stackrel{(\ref{eq:com_inv})}{=}(K^{-1})_{ij}$.
\end{proof}}

\cor{\label{cor:invK}For $i,j\in {\cal Y}$ ($i\neq j$) $\bY_i$ and $\bY_{j}$ are independent 
given ${\cal Y}\setminus \{i,j\}\subset \bY$ if and only if 
$$K_{ij}^{-1}=0.$$ }\rm

\noindent \emph{Remark.}  Kulesza in \cite{kul} found that $i\in \bY \perp j\in \bY| ({\cal Y}\setminus\{i,j\}\subset \bY)$
if and only if $K_{ij}^{-1}=0$.\vspace{0.2cm}

By  Proposition \ref{prop:sumrez} (c) ${\cal Y}\setminus \bY$ is a DPP with the kernel $I-K$. 
But the more interesting thing is that ${\cal Y}\setminus \bY$ is the $L$-ensemble with the kernel
\begin{equation}
 \bar{L}=K^{-1}-I.
\end{equation}

Now, the Corollary \ref{cor:invK} can be restated in the terms of the matrix $\bar{L}$.

\cor{For $i,j\in {\cal Y}$ ($i\neq j$) $\bY_i$ and $\bY_{j}$ are independent 
given ${\cal Y}\setminus \{i,j\}\subset \bY$ if and only if 
$$\bar{L}_{ij}=0.$$} \rm

Looking at the process $\bY={\cal Y}\setminus ({\cal Y}\setminus \bY)$ we have

\cor{For $i,j\in {\cal Y}$ ($i\neq j$) $\bY_i$ and $\bY_{j}$ are independent 
given $({\cal Y}\setminus \{i,j\})\cap  \bY=\emptyset$ if and only if 
$$L_{ij}=0.$$} \rm

\section{Comparison to Gaussian graphical models}
The way independence is encoded in matrices $K$ and $L$ is similar 
to way independence is encoded in covariance matrix $\Sigma$ and 
precision matrix $\Sigma^{-1}$ of the Gaussian random vector.\vspace{0.2cm} 

The question is, can we, from the structure of the matrix $L$, say 
more about conditional independencies in a DPP? Is there a similar
result as in the Gaussian graphical models?\vspace{0.3cm}

We will briefly review the results we have in Gaussian graphical models.

We will assume $V=\{1, \ldots , n\}$ and let the process $$X=(X_v:v\in V)$$ be a 
a normal random vector with expectation $\mu$ and a positive definite covariance matrix $\Sigma$.

\defi{For a symmetric matrix $M$ we will say that $G_M=(V,E_{M})$ is a graph 
induced by the matrix $M$ if the set of edges is given by
$$E_M=\{\{i,j\}\,:\, M_{ij}\neq 0, i\neq j\}.$$ }\rm

The following results are well known for Gaussian random vectors.

\thm{\label{thm:gss}\begin{enumerate}[(a)]
    \item For disjoint subsets  $A,B,C$ of $V$ 
    $$X_A\perp X_B|X_C$$
if and only if $\Sigma^C_{AB}={\bf 0}$.
      \item \label{thm:gss:b}For $k,j\in V$ with $k\neq j$ 
$$X_k\perp X_j|X_{V\setminus\{k,j\}}$$
if and only if $\Sigma^{-1}_{k,j}=0$.
     \end{enumerate}
}


\defi{\begin{enumerate}[(a)]
       \item We say that the process $X$ has the pairwise Markov property with 
respect to the structure of the graph $G=(V,E)$ if 
$X_k\perp X_j|X_{V\setminus\{k,j\}}$ holds for all $\{k,j\}\notin E$.
      \item We say that the process $X$ has the global Markov property if
for $A,B,C$ are disjoint subsets of $V$ such that $C$ separates $A$ and $B$,
i.e. any path starting at a vertex in $A$ and ending in $B$ has to go through 
a vertex in $C$, we have $X_A\perp X_B|X_C$.
      \end{enumerate}
}\rm \vspace{0.2cm}

The following is a consequence of the famous Hammeresley-Clifford Theorem
and the fact that $X$ has a positive density. (See \S3.2.1. and Theorem 3.9. in \cite{grph:model}.)
\thm{\label{glb-prw}The process $X$ has the pairwise Markov property with respect
to graph $G=(V,E)$ if and only if it has the global 
Markov property with respect to $G$.}

\cor{$X$ is a has the pairwise 
Markov property with the respect to the structure of the graph $G_{\Sigma^{-1}}=(V,E_{\Sigma^{-1}})$.
Further, $X$ also has the global Markov property with the respect to $G_{\Sigma^{-1}}$.
\begin{proof}
 From the definition, using Theorem \ref{thm:gss}. (\ref{thm:gss:b}) the pairwise property 
follows. The global property follows from Theorem \ref{glb-prw}. 
\end{proof}
}

\thm{\label{thm:Mstr}Let $M$ be a positive definite $n\times n$ matrix, and $G_{M^{-1}}=(V, E_{M^{-1}})$
a graph induced by $M^{-1}$. If $A,B,C$ are disjoint subsets of $V$ such that $C$ separates $A$ and $B$, then 
$$M^{C}_{AB}={\bf 0}.$$
\begin{proof}
Let $Y\sim N(0,M)$. 
By Theorem \ref{glb-prw}, $Y$ has the global Markov property with respect to the graph $G_{M^{-1}}$. 
Hence $Y_A$ is independent of $Y_B$ given $Y_C$, and by Theorem \ref{thm:gss}.(a) this is true if and only if 
$M^{C}_{AB}={\bf 0}$.
\end{proof}}\rm

\section{Graphs induced by the $L$-ensemble}

From the structure of the $L$-ensemble we can get some information
about other conditional independencies. The following is a version 
of the global Markov property for $L$-ensembles.

\thm{\label{thm:Lgrph}Let the determinantal point process $\bY$ be an $L$-ensemble and $G_L$ be a graph induced by the kernel $L$. 
If $A,B,C$ are disjoint subsets of $V$ such that $C$ separates $A$ and $B$,
then $\bY_A$ is independent of $\bY_B$ given that $\bY \cap C=\emptyset$.

\begin{proof}
$L$ has off-diagonal zeros in the same places as $(I-K)^{-1}$ (see $(\ref{eq:L})$). By Theorem \ref{thm:Mstr}, we have 
that $(I-K)^C_{AB}=0$. Hence, by Theorem \ref{thm:mainstrong}, $({\cal Y}\setminus \bY)\cap A$ and $({\cal Y}\setminus \bY)\cap B$ 
are independent given $C\subset {\cal Y}\setminus \bY$. Hence, the claim follows. 
\end{proof}
}

\thm{\label{cor:Lgrph2}Let the determinantal point process $\bY$ be an $L$-ensemble and $G_L$ be a graph induced by the kernel $L$. 
If $A,B,C,D$ are disjoint subsets of $V$ such that $C$ separates $A$ and $B$,
then $\bY_A$ is independent of $\bY_B$ given that $\bY \cap C=\emptyset$ and $D\subset \bY$.
%
\begin{proof}
Let $D_A$ be all vertices from $D$ that are connected to vertices from $A$ with paths in 
$G_L$ that do not pass through $C$. We set $D'_A=D\setminus D_A$.
For $S_A\in \sigma (\bY_A)$ and $S_B\in \sigma (\bY_B)$, using 
the fact that $C$ separates $A\cup D_A$ and  $B\cup D_A'$ and Theorem \ref{thm:Lgrph}.
we have
\begin{align*}
 & \Pr(S_A\cap S_B | \bY \cap C=\emptyset,D\subset \bY)\\
 =& \Pr(S_A\cap S_B, D\subset \bY | \bY \cap C=\emptyset)/\Pr(D\subset \bY | \bY \cap C=\emptyset)\\
 =& \frac{\Pr(S_A, D_A\subset \bY | \bY \cap C=\emptyset)\Pr(S_B, D'_A\subset \bY | \bY \cap C=\emptyset)}{\Pr(D\subset \bY | \bY \cap C=\emptyset)}\\
=& \frac{\Pr(S_A, D_A\subset \bY | \bY \cap C=\emptyset)\Pr(D'_A\subset \bY | \bY \cap C=\emptyset)}{\Pr(D'_A\subset \bY | \bY \cap C=\emptyset)\Pr(D_A'\subset \bY | \bY \cap C=\emptyset)}\\
 & \cdot \frac{\Pr(S_B, D'_A\subset \bY | \bY \cap C=\emptyset)\Pr(D'_A\subset \bY | \bY \cap C=\emptyset)}{\Pr(D\subset \bY | \bY \cap C=\emptyset)}\\
 =& \frac{\Pr(S_A, D\subset \bY | \bY \cap C=\emptyset)}{\Pr(D\subset \bY | \bY \cap C=\emptyset)}\cdot \frac{\Pr(S_B, D\subset \bY | \bY \cap C=\emptyset)}{\Pr(D\subset \bY | \bY \cap C=\emptyset)}\\
& = \Pr(S_A|\bY \cap C=\emptyset,D\subset \bY)\Pr(S_B|\bY \cap C=\emptyset,D\subset \bY).
\end{align*}

\end{proof}
}\rm

\prop{\label{prop:Lgrph}Let the determinantal point process $\bY$ be an $L$-ensemble and $G_L$ be a graph induced by the kernel $L$. 
Let\begin{itemize}
   \item $A_1$, \ldots $A_n$, $C$ and $D$ are disjoint subsets of ${\cal Y}$;
   \item $C$ separates sets $A_1$, \ldots $A_n$ in $G_L$.
  \end{itemize}
Then $\bY_{A_1}$, \ldots, $\bY_{A_n}$ are independent given $C\cap \bY=\emptyset$
and $D\subset \bY$.}\rm

\section{Final remarks}
Proposition \ref{prop:main} gives necessary and sufficient conditions 
for conditional independencies, but it is not easy to practically check them.
Further, estimating $K$ is conjectured to be an NP-hard problem (\cite{kul_task}).\vspace{0.2cm} 

On the other hand, Theorem \ref{thm:Lgrph} gives us only sufficient conditions on the kernel $L$   and given a sparse matrix $L$ we can read many conditional independencies
from its structure without any additional matrix transformations. Further, there are ways to estimate kernel $L$ (\cite{kul_task}).\vspace{0.2cm} 

Although the independence induced by the graph structure is not as
strong as in the case of graphical models, it still provides important 
information about the process and is useful for better understanding 
of this process.


\end{document}